\documentclass[11pt]{amsart}
\usepackage{amssymb,latexsym,comment,url}
\usepackage{amsfonts}
\usepackage{enumerate}

\usepackage[mathscr]{eucal}
\usepackage{xcolor}
\usepackage{eqlist}

\newcommand{\Orb}{\operatorname{Orb}}

\theoremstyle{plain}

\theoremstyle{definition}

\theoremstyle{remark}

\numberwithin{equation}{section}

\newcommand{\PGL}{\operatorname{PGL}}
\newcommand{\Aut}{\operatorname{Aut}}
\newcommand{\Char}{\operatorname{char}}

\newcommand{\Q}{\mathbb{Q}}

\newcommand{\PP}{\mathbb{P}}
\newcommand{\Magma}{{\sf MAGMA }}

\newenvironment{Proof}{\par\noindent{\sc Proof:}}%
                      {\hspace*{\fill}\nobreak$\Box$\par\medskip}
\newenvironment{ProofOf}[1]{\par\noindent{\sc Proof of #1:}}%
                       {\hspace*{\fill}\nobreak$\Box$\par\medskip}

\newtheorem{Proposition}{Proposition}[section]
\newtheorem{Theorem}[Proposition]{Theorem}
\newtheorem{Lemma}[Proposition]{Lemma}

\newtheorem{Example}[Proposition]{Example}
\newtheorem{Conjecture}[Proposition]{Conjecture}
\newtheorem{Definition}[Proposition]{Definition}
\newtheorem{Remark}[Proposition]{Remark}

\begin{document}

   \title[Simultaneous periodic points of quadratic maps]{Simultaneous Rational Periodic Points of Degree-2 Rational Maps}

\author[Burcu Barsak\c{c}i]%
{Burcu Barsak\c{c}i}
\address{Faculty of Engineering and Natural Sciences, Sabanc{\i} University, Tuzla, \.{I}stanbul, 34956 Turkey}
\email{burcubarsakci@sabanciuniv.edu}

\author[M. Sadek]%
{Mohammad~Sadek}
\address{Faculty of Engineering and Natural Sciences, Sabanc{\i} University, Tuzla, \.{I}stanbul, 34956 Turkey}
\email{mohammad.sadek@sabanciuniv.edu}

\date{\today}

\maketitle
	
\begin{abstract}
Let $S$ be the collection of quadratic polynomial maps, and degree $2$-rational maps whose automorphism groups are isomorphic to $C_2$ defined over the rational field. Assuming standard conjectures of Poonen and Manes on the period length of a periodic point under the action of a map in $S$, we give a complete description of triples $(f_1,f_2,p)$ such that $p$ is a rational periodic point for both $f_i\in S$, $i=1,2$.  We also show that no more than three quadratic polynomial maps can possess a common periodic point over the rational field.
In addition, under these hypotheses we show that two nonzero rational numbers $a, b $ are periodic points of the map $\phi_{t_1,t_2}(z)=t_1  z + t_2/z$ 
for infinitely many nonzero rational pairs $(t_1, t_2)$ if and only if $a^2 = b^2$. 
\end{abstract}

\section{Introduction}

A (discrete) dynamical system consists of a set $S$ and a map $\phi: S \rightarrow S$ that permits iteration
$$ \phi^n =\underbrace{ \phi \circ \phi \circ \ldots \circ \phi }_{n\textrm{- times}}.$$
For a given point $P \in S$, the orbit of $P$ is the set
$$ \mathcal{O}_{\phi} (P)= \mathcal{O} (P)= \{ \phi^n (P): n \geqslant 0 \}.$$
One of the main goals of dynamics is to classify the points $P$ in the set $S$ according to the size of their orbits $\mathcal{O}_{\phi} (P)$. If $\mathcal{O}_{\phi} (P)$ is infinite, $P$ is called a wandering point; otherwise $P$ is called a preperiodic point. In this paper, we focus on a special type of preperiodic points, the so called periodic points. A point $P \in S$ is periodic if there exists an integer $n > 0$ such that $\phi^n (P)= P$, where $n$ is called the period of $P$. If $n$ is the smallest such integer, we say that $P$ has exact period $n$. In our work, we will take the set $S$ to be the rational field $\mathbb{Q}$ and the map $\phi$ to be a degree-$2$ rational map of either forms
$$f(z)=z^2 + c, \qquad   \phi_{k, b}(z)=kz + \frac{b}{z},$$
where $k, b, c \in \mathbb{Q}$ with $k, b$ are nonzero.
Due to the work of Northcott, \cite{Northcott}, we know that these maps can have only finitely many rational periodic points.

A complete classification of rational quadratic polynomials $f(z)$ with periodic points of exact period $1, 2$ or $3$ can be found for example in \cite{Poonen, WaldeRusso}. A polynomial $f(z)$ cannot have a rational periodic point of exact period $4$, \cite{Morton}, or exact period $5$, \cite{fps}. Under the assumption of Birch-Swinnerton-Dyer conjecture, Stoll proved that a polynomial $f(z)=z^2+c$ cannot have a rational periodic point of exact period $6$, \cite{Stoll}. Poonen  conjectured that a rational quadratic polynomial $f(z)$ cannot have a rational periodic point of exact period $\geqslant 6$, \cite{Poonen}. In this work, we assume that the latter conjecture holds. 

In \cite{Manes}, a classification of rational maps of the form $\phi_{k,b}(z)$ that possess rational periodic points of exact period $1, 2$ or $4$ is given. It was proved that $\phi_{k,b} (z)$ cannot  have a rational periodic point of exact period $3$. It was also conjectured that $\phi_{k,b}(z)$ cannot have a rational periodic point of exact period at least $5$. Again, in this work we assume that the latter conjecture holds.

We classify all triples of rational numbers $(k,b,c)$ such that $f(z)=z^2 + c$ and $ \phi_{k, b}(z)=kz + \frac{b}{z}$ possess a rational periodic point $p$ under the assumption of Poonen and Manes' conjectures. In addition, we find all such triples for which
$$ | \Orb_{f} (p) \cap \Orb_{ \phi_{k, b}} (p) | \geqslant 2.$$
 We also prove that the size of the intersection of $\Orb_{f} (p)$ and $\Orb_{ \phi_{k, b}} (p)$ cannot exceed $2$. 
 
Under the assumption of Manes' conjecture, we list down the tuples $(k_1, b_1, k_2, b_2)$ such that
the maps $\phi_{k_1, b_1}(z)$ and $\phi_{k_2, b_2}(z)$ have a simultaneous rational periodic point $p$. We then illustrate which such tuples allow the following
$$ | \Orb_{ \phi_{k_1, b_1}} (p) \cap \Orb_{ \phi_{k_2, b_2}} (p) | \geqslant 2.$$
 Again, we show that the above intersection size cannot exceed $2$ unless $\phi_{k_1, b_1} = \pm \phi_{k_2, b_2}$.

In \cite{BakerDeMarco}, it was proved that for fixed $c_1, c_2 \in \mathbb{C}$, the set of $t \in \mathbb{C}$ such that both $c_1$ and $c_2$ are preperiodic for $z^d + t$ is infinite if and only if $c_1 ^d = c_2 ^d$ where $d$ is greater than one. Assuming Manes' conjecture, we prove a similar result over the rational field. 
Setting
$\phi_{t_1, t_2}(z)=t_1  z + \frac{t_2}{z},$
two rational numbers $a, b \in \mathbb{Q}^{\times}$ are periodic points of $\phi_{t_1, t_2}(z)$ for infinitely many nonzero rational pairs $(t_1, t_2) $ if and only if $a^2 = b^2$.

Although it is easily seen that given a nonzero rational number $p$ we can find infinitely many nonzero rational pairs $(k,b)$ such that $p$ is a periodic point of $\phi_{k, b}(z)$, there are at most  three quadratic polynomial maps sharing $p$ as a common periodic point.  The latter result is in the same realm of results on the size of sets $S$ of quadratic polynomial maps for which $p$ is preperiodic for any composition of maps in $S$, \cite{Hindes}.

\subsection*{Acknowledgments} This work is supported by The Scientific and Technological Research Council of Turkey, T\"{U}B\.{I}TAK; research grant: ARDEB 1001/120F308 and BAGEP Award of the Science Academy, Turkey.

\section{Periodic points of degree-$2$ rational maps}

In this section, we collect some preliminaries on dynamical systems and known results on periodic points of the quadratic rational maps that we deal with in this article.  
\subsection{Preliminaries}
We start with recalling the definition of the dynatomic polynomial associated to a certain rational map defined over a field $K$ with algebraic closure $\overline{K}$.

\begin{Definition} \cite{dynamicalbook} \label{Dynatomic}
Let $K$ be a field. Let $\phi(z) \in K(z)$ be a rational function of degree $d$. We write
$$ \phi^n =[F_n (x, y), G_n (x, y)] $$ 
where $F_n, G_n \in K[x, y]$ are homogeneous polynomials of degree $d^n$ for any $n\geqslant0$. The $n$-period polynomial of $\phi$ is the polynomial 
$$ \Phi_{\phi, n} (x, y)=y \cdot F_n (x, y) - x \cdot G_n (x, y). $$\\
The $n^{th}$ dynatomic polynomial of $\phi$ is the polynomial
$$ \Phi^*_{\phi, n} (x, y)=\prod_{k|n} (y \cdot F_k (x, y) - x \cdot G_k (x, y))^{\mu(n/k)}=\prod_{k|n} (\Phi_{\phi, k} (x, y))^{\mu(n/k)} $$
where $\mu$ is the M\"{o}bius function, defined by $ \mu(1)=1$, $\mu(n)=(-1)^l$ if $n=p_1 p_2 \ldots p_l$ with $p_i$ distinct primes, and $\mu(n)=0$ if $n$ is not square free.\\
If $\phi$ is fixed, we write $\Phi_n$ and $\Phi^*_n$. We will write $\Phi_n (z)$ for $\Phi_n(z,1)$ and $\Phi^*_n (z)$ for $\Phi^*_n(z,1)$.     
\end{Definition}
By definition, every exact period-$n$ point of $\phi$ is a root of $\Phi^*_{\phi, n} (x,y)$. It is important to notice that $\Phi^*_{\phi, n} (x,y)$ can have roots whose periods divide $n$ and strictly smaller than $n$. 

Let $f \in \PGL_2(\overline K)$ act on the points of $\PP^1$ as a fractional linear transformation in the usual way. Then we define the rational map
$\phi^f = f^{-1}\phi f$. In fact, two rational maps $\phi$ and $\psi$ are said to be {\em linearly conjugate} if there is
some $f\in  \PGL_2(\overline K)$ such that $\phi^f = \psi$. They are linearly conjugate over $K$ if there is some 
$f \in \PGL_2(K)$ such that $\phi^f = \psi$.

We notice that if $P$ is a periodic point of exact period $n$ for $\phi$, then $f^{-1}(P)$ is a periodic point of exact period $n$ for $\phi^f$. One can argue similarly for preperiodic points of $\phi$ and $\phi^f$. In other words, linear conjugation preserves the dynamical behavior of points under rational maps. Moreover, if $\phi, f$ and $P$ are defined over $K$ such that $\phi^n(P)=P$, then $\psi:=\phi^f$ and $Q:=f^{-1}(P)$ are defined over $K$ with $\psi^n(Q)=Q$.

\subsection{Quadratic polynomial maps} 

 Any quadratic polynomial map $\phi(z)=Az^2+Bz+C\in K[z]$,  $A\in K^{\times}$, is linearly conjugate over $K$ to a map of the form $\psi(z)=z^2+c$ for some $c\in K$. Therefore, we only consider quadratic polynomial maps of this form. 

If $K$ is chosen to be the rational field $\Q$,  a complete classification of quadratic polynomial maps with periodic points of periods $1,2,$ or $3$ was given in \cite{WaldeRusso}. The following can be found for example as \cite[Theorem 1]{Poonen}.
\begin{Proposition}
\label{prop:1}
Let $f(z)=z^2 +c$ with $c\in \Q$. Then
\begin{itemize}
\item[1)] $f(z)$ has a rational point of period $1$, i.e., a rational fixed point, if and only if $c = 1/4 -\rho^2$ for some $\rho\in  \Q$. In this case, there are exactly two, $1/2 +\rho$ and $1/2 -\rho$, unless $\rho = 0$, in which case they coincide.
\item[2)] $f(z)$ has a rational point of period $2$ if and only if $c = -3/4 -\sigma^2$  for some $\sigma\in\Q$, $\sigma\ne 0$. In this case, there are exactly two, $-1/2 +\sigma$ and $-1/2 -\sigma$  (and these form a $2$-cycle).
\item[3)] $f(z)$ has a rational point of period $3$ if and only if
$$c = -\frac{\tau^6  + 2\tau^ 5 + 4\tau^4 + 8\tau^3 + 9\tau^2 + 4\tau + 1}{
4\tau^2(\tau +1)^2}$$
for some $\tau  \in \Q$, $\tau\ne  -1, 0$. In this case, there are exactly three,
\begin{eqnarray*}
x_1 = \frac{\tau^3 +2\tau^2 +\tau +1}{ 2\tau(\tau +1)},\quad
x_2 = \frac{\tau^3 -\tau -1}{ 2\tau(\tau +1)},\quad
x_3 = -\frac{ \tau^ 3 + 2\tau^2 + 3\tau + 1}{ 2\tau(\tau +1)}
\end{eqnarray*}
and these are cyclically permuted by $f(z)$.
\end{itemize}
\end{Proposition}   
The following conjecture can be found in \cite{fps,Poonen}.
\begin{Conjecture}\label{Con:Poonen}
If $N \ge 4$, then there is no quadratic polynomial $f (z) \in\Q[z]$ with a rational point of exact period $N$ .
\end{Conjecture}
The conjecture has been proved for $N=4$, \cite{Morton}, for $N=5$, \cite{fps}, and conditionally for $N=6$, \cite{Stoll}. 
\subsection{Degree-$2$ rational maps}
Given a rational map $\phi(z)\in K(z)$, we define the automorphism group of $\phi$, $\Aut(\phi)$, to be $\{f\in\PGL_2(\overline K):\phi^f=\phi\}$. Let $\phi$ be a degree-$2$ rational map defined over $K$, $\Char K\ne 2,3$. Then $\Aut(\phi)\cong C_2$ if and only if $\phi$ is linearly conjugate over $K$ to a  map of the form $\phi_{k,b}(z)=kz+b/z$, $k\not\in\{0,-1/2\}$, $b\in K^{\times}$. In addition, two such maps $\phi_{k,b}$ and $\phi_{k',b'}$ are linearly conjugate over $K$ if and only if $k=k'$ and $b/b'\in (K^{\times})^2$. The map $\phi_{k,b}$ has the automorphism $z\mapsto -z$, see \cite[Lemma 1]{Manes}. Since the homogeneous form of $\phi_{k,b}(z)=kz+b/z$ is given by $[kx^2+by^2:xy]$, it follows that the point at infinity is a $K$-rational fixed point for $\phi_{k,b}(z)$. Therefore, throughout the sequel when we mention periodic points of $\phi_{k,b}(z)$ we mean finite periodic points. 

The following proposition describes $\Q$-rational periodic points of maps of the form $\phi_{k,b}(z)=kz+b/z$, $k,b\in\Q^{\times}$, see \cite[\S 2]{Manes}.     

\begin{Proposition}
\label{prop:2}
Let $\phi_{k,b}(z)=kz+b/z$, $k\in \Q^{\times}$, $b\in \Q^{\times}/(\Q^{\times})^2$.
\begin{itemize}
\item[1)] $\phi_{k,b}(z)$ has a $\Q$-rational fixed point if and only if $b/(1-k)=m^2$ for some $m\in \Q^{\times}$, $k\ne 1$. In this case there are two finite $\Q$-rational fixed points, $\pm m$.
\item[2)] $\phi_{k,b}(z)$ has a $\Q$-rational periodic point of exact period $2$ if and only if $b/(k+1)=-m^2$ for some $m\in \Q^{\times}$, $k\ne -1$.  In this case there are exactly two such points, $\pm m$.
\item[3)] $\phi_{k,b}(z)$ has no $\Q$-rational periodic point of exact period $3$.
\item[4)] $\phi_{k,b}(z)$  has a $\Q$-rational periodic point of exact period $4$ if and only if $$k=\frac{2m}{m^2-1},\qquad b=\frac{-m}{m^4-1}\qquad\textrm{for some }m\in\Q\setminus\{0,\pm 1\}.$$ In this case, the rational periodic points are $$x_1=\frac{1}{m^2+1},\quad x_2=\frac{-m}{m^2+1}, \quad x_3=\frac{-1}{m^2+1}, \quad x_4=\frac{m}{m^2+1}.$$
\end{itemize}
\end{Proposition}  
The following conjecture is \cite[Conjecture 1]{Manes}.
\begin{Conjecture}\label{Con:Manes}
If $\phi(z) \in\Q(z)$ is a degree-$2$ rational map with $\Aut(\phi)\cong  C_2$, then $\phi$ has
no rational point of exact period $N > 4$.
\end{Conjecture}

\section{Simultaneous rational periodic points of $f(z)$ and $\phi_{k,b}(z)$}

\begin{Proposition} \label{mixed1}
 Let $f(z)=z^2+c$ and $\phi_{k,b}(z)=kz+b/z$, $c\in\Q, k,b\in\Q^{\times}$. Given $p\in\Q^{\times}$, the set of triples $(k,b,c)$ such that $p$ is a rational periodic point of exact period $1$ (rational fixed point) of $f(z)$ and a rational periodic point of exact period $n$ of $\phi_{k,b}(z)$ is described as follows 
\[
\begin{array}{cc}
      (k,b,c)=\left(\frac{ q+p}{p},-qp,p-p^2\right),\quad  q \in \mathbb{Q}\setminus\{0,-p\}& \textrm{ if } n=1 \\
      (k,b,c)=\left(\frac{ q-p}{p},-qp,p - p^2\right),\quad
 q \in \mathbb{Q}\setminus\{0,p\}& \textrm{ if } n=2 \\
 (k,b,c)=\left(\frac{2m}{m^2 -1},- \frac{p^2  (m^2 + 1)}{m (m^2 -1)}, p - p^2\right), \quad m \in \mathbb{Q} \setminus \{0, \pm1\} & \textrm{ if } n=4 \\
    \end{array}
\]
\end{Proposition}
\begin{Proof}
One may check easily that $p$ is a rational fixed point for both $f(z)$ and $\phi_{k,b}(z)$ if $k, b, c$ are as given in the statement of the proposition. Now,  we know that $f(z)$ has a rational fixed point if and only if $c=1/4 - \rho^2$ for some $\rho \in \mathbb{Q}$, see Proposition \ref{prop:1}, where $1/2\pm \rho$ are rational fixed points of $f(z)$. Suppose that $p=1/2 + \rho$ is a rational fixed point of $\phi_{k,b}(z)$, it follows that
$
1/2+\rho=\frac{ \pm  \sqrt{b -bk}}{ -1 + k},
$
where $k \neq 1$ and
$
b-bk = q^2 
$
for some $q \in \mathbb{Q}^{\times}$.  It follows that
$
k=\frac{\pm 2q+2\rho+1}{2\rho +1}
$ and
$
b=\mp (1+2\rho)q/2.
$
We obtain the same parametrization if we take $p=1/2 - \rho$ to be the common periodic point of $f(z)$ and $\phi_{k,b}(z)$.

For $n=2$, the proof is similar. 

As for $n=4$,  we recall that the $4^{th}$-dynatomic polynomial of $\phi_{k,b}(z)$ factors as follows $$
\Phi^*_4 (z)=\Psi_4 (z) \cdot \Lambda_4 (z) 
$$ where 
\begin{eqnarray*}
\Psi_4 (z)= b^2 k &+& 2 b z^2 + 2 b k^2 z^2 + k z^4 + k^3 z^4,\\
  \Lambda_4 (z)= b^4 k^5 &+& b^3 z^2 + b^3 k^2 z^2 + 2 b^3 k^4 z^2 + 4 b^3 k^6 z^2 + 
  b^2 k z^4 + 3 b^2 k^3 z^4 + 4 b^2 k^5 z^4 + 6 b^2 k^7 z^4\\ &+& 
  b k^4 z^6 + 2 b k^6 z^6 + 4 b k^8 z^6 + k^9 z^8.
  \end{eqnarray*} Moreover,  $\Lambda_4 (z)$ does not have a rational root, see \cite[\S 2]{Manes}. If $1/2  +  \rho$ is a rational periodic point of exact period $4$ of $\phi_{k,b}(z)$, then it is a root of $\Psi_4 (z)$, it follows that $\rho$ is one of the following four expressions 
\begin{equation*} \label{rho1}
\frac{1}{2} \left( -1 \pm \sqrt{ \frac{-4b \cdot (1 + k^2   \pm \sqrt{1 + k^2})}{k \cdot (1+k^2)}}  \right).
\end{equation*}

The rationality of $\rho$ implies that $1+k^2$ is a rational square, hence $k=\frac{2m}{m^2-1}$ for some $m\in\Q\setminus\{0,\pm 1\}$.  Moreover, $$b=- \frac{s^2 \cdot (m^2 + 1)}{4m (m^2 -1)},\qquad\textrm{ where }s= \sqrt{ \frac{-4b \cdot (1 + k^2   \pm\sqrt{1 + k^2})}{k \cdot (1+k^2)}}.$$
Now, the expressions for $b$ and $c$ follow by noticing that $p=1/2+\rho=\pm s/2$.
In fact, the cycle of periodic points of period $4$ of $\phi_{k,b}(z)$ is given by
$ \left(p,p/m, -p, -p/m \right)$.
\end{Proof}
\begin{Example}
The point $p=3/2$ is a common fixed point of 
$ f(z)=z^2 -\frac{3}{4},$
and $\phi(z)=\frac{5z}{3}-\frac{3}{2z}$ corresponding to $q=1$ in Proposition \ref{mixed1}.

Setting $p=3$ and $q=\frac{1}{2}$, we get that $3$ is a rational fixed point of $ f(z)=z^2 - 6, $ whereas it is a periodic point of exact period $2$ of $\phi(z)=-\frac{5z}{6}-\frac{3}{2z}.$ 

Finally, taking $m=2$ and $p=2$ in Proposition \ref{mixed1} implies that $2$ is a rational fixed point of 
$ f(z)=z^2 - 2, $ while it is a periodic point of exact period $4$ for
$\phi(z)=\frac{4z}{3}-\frac{10}{3z}.$
\end{Example}

\begin{Proposition} \label{mixed2}
Let $f(z)=z^2+c$ and $\phi_{k,b}(z)=kz+b/z$, $c\in\Q, k,b\in\Q^{\times}$. Given $p\in\Q\setminus\{0,-1/2\}$, the set of triples $(k,b,c)$ such that $p$ is a rational periodic point of exact period $2$ of $f(z)$ and a rational periodic point of exact period $n$ of $\phi_{k,b}(z)$ is described as follows 
\[
\begin{array}{cc}
      (k,b,c)=\left(\frac{ q+p}{p},-qp,-(p^2+p+1)\right),\quad  q \in \mathbb{Q}\setminus\{0,-p\}& \textrm{ if } n=1 \\
      (k,b,c)=\left(\frac{ q-p}{p},-qp,-(p^2 +p+1)\right),\quad
 q \in \mathbb{Q}\setminus\{0,p\}& \textrm{ if } n=2 \\
 (k,b,c)=\left(\frac{2m}{m^2 -1},- \frac{p^2  (m^2 + 1)}{m (m^2 -1)}, -(p^2+p+1) \right), \quad m \in \mathbb{Q} \setminus \{0, \pm1\} & \textrm{ if } n=4 \\
    \end{array}
\]
\end{Proposition}
\begin{Proof}
The proof is similar to the proof of Proposition \ref{mixed1} using Propositions \ref{prop:1} and \ref{prop:2}.
\end{Proof}
\begin{Example}
The point $p=1/2$ is a periodic point of exact period $2$ of 
$ f(z)=z^2 -\frac{7}{4},$
and a fixed point of $\phi(z)=3z-\frac{1}{2z}$ corresponding to $q=1$ in Proposition \ref{mixed2}.

Setting $p=1$ and $q=-1$, we get that $1$ is a periodic point of exact period $2$ of $ f(z)=z^2 - 3 $ and $\phi(z)=-2z+\frac{1}{z}.$ 

Finally, taking $m=3$ and $p=-1$ in Proposition \ref{mixed2} implies that $-1$ is a periodic point of exact period $2$ of 
$ f(z)=z^2 - 1, $ while it is a periodic point of exact period $4$ for
$\phi(z)=\frac{3z}{4}-\frac{5}{12z}.$
\end{Example}

\begin{Proposition} \label{mixed3}
Let $f(z)=z^2+c$ and $\phi_{k,b}(z)=kz+b/z$, $c\in\Q, k,b\in\Q^{\times}$.  The set of triples $(k,b,c)$ such that $f(z)$ has a rational periodic point $x_i$ of exact period $3$ which is a rational periodic point of exact period $n$ of $\phi_{k,b}(z)$ is described as follows 
\[
\begin{array}{cc}
      (k,b,c)=\left(1-q,q  x_i^2,- \frac{\tau^6 + 2 \tau^5 + 4 \tau^4 + 8 \tau^3 + 9 \tau^2 + 4 \tau +1}{4 \tau^2 (\tau +1)^2}\right)& \textrm{ if } n=1 \\
      (k,b,c)=\left(q-1,-qx_i^2,- \frac{\tau^6 + 2 \tau^5 + 4 \tau^4 + 8 \tau^3 + 9 \tau^2 + 4 \tau +1}{4 \tau^2 (\tau +1)^2}\right)& \textrm{ if } n=2 \\
       \end{array}      
\]
where $q \in \mathbb{Q} \setminus \{ 0, 1\}$ and $\tau \in \mathbb{Q} \setminus \{ -1, 0\}$, $i=1,2,3$,
\[
\begin{array}{cc}
 (k,b,c)=\left(\frac{2m}{m^2 -1},- x_i^2  \frac{(m^2 +1)}{m (m^2 -1)},- \frac{\tau^6 + 2 \tau^5 + 4 \tau^4 + 8 \tau^3 + 9 \tau^2 + 4 \tau +1}{4 \tau^2 (\tau +1)^2}\right), \quad m \in \mathbb{Q} \setminus \{0, \pm1\}, i=1,2,3 & \textrm{ if } n=4 \\
    \end{array}
\]
where \[x_1=\frac{\tau^3 + 2 \tau^2 + \tau + 1}{2 \tau (\tau +1)},\,
x_2=\frac{\tau^3 - \tau - 1}{2 \tau (\tau +1)},\,
x_3=-\frac{\tau^3 + 2 \tau^2 + 3 \tau + 1}{2 \tau (\tau +1)},\qquad\textrm{ where }
\tau\in\Q\setminus\{0,-1\}.\]
\end{Proposition}
\begin{Proof}
The cases $n=1,2$ can be treated similarly as in Propositions \ref{mixed1} and \ref{mixed2}. 

For $n=4, $ Proposition \ref{prop:1} asserts that $f(z)$ has a rational periodic point of exact period $3$ if and only if $c$ and $x_i$ are given in the statement of the proposition. To force one of these points to be a periodic point of exact period $4$ for $\phi_{k,b}$ for some rational numbers $k,b$, the point should be a root of $\Psi_4 (z)$, see the proof of Proposition \ref{mixed1}. The four roots of $\Psi_4 (z)$ are given by 
\begin{equation}
\nonumber
 \pm \sqrt{\frac{-b}{k} \left( 1\pm \frac{1}{\sqrt{1 +k^2}}\right)}. 
\end{equation}
The rationality of the periodic points yields the result by proceeding as in the proof of Proposition \ref{mixed1}.  
Notice that in this case the cycle of periodic points of exact period $4$ of $\phi_{k,b}(z)$ is given by 
$ \left( x_i,  x_i/m, - x_i, - x_i/m \right).$
\end{Proof}
\begin{Example}
Fixing $\tau=1$ and $q=16$, then $x_2=-1/4$ is a rational periodic point of exact period $3$ of $f(z)=z^2 -29/16$ and it is a rational fixed point of
$\phi(z)=-15z + \frac{1}{z}$. 

Setting $\tau=1/2$ and $q=9$, then $x_1=17/12$ is a rational periodic point of exact period $3$ of $f(z)=z^2 -421/144$ while it is  
 a rational periodic point of exact period $2$ of
$\phi(z)=8z - \frac{289}{16z}$.

If we set $\tau=-1/2$ and $m=2$,  then $x_3=-1/4$ is a rational periodic point of exact period $3$ of $f(z)=z^2 -29/16$ and it is a rational point of exact period $4$ of
$\phi(z)=\frac{4z}{3} - \frac{5}{96z}$.
\end{Example}
The propositions above imply the following result. 
\begin{Theorem}
\label{Thm1}
Assuming Conjecture \ref{Con:Poonen} and Conjecture \ref{Con:Manes}, the following table contains all pairs $f(z)=z^2+c$ and $\phi_{k,b}(z)=kz+b/z$, $c\in\Q$, $k,b \in\Q^{\times}$, such that some rational number $p$ is a periodic point for both maps.
 \begin{center}
\begin{tabular}{ | m{9em} | m{3.5cm}| m{0.3cm}| m{2.7cm}| m{1.7cm}| m{0.3cm}|  } 
\hline
$\phi_{k, b}(z)$ & Periodic Points & PL  &$\beta(z)$ & Periodic Points & PL \\ 
\hline
$\frac{q+p}{p} z -\frac{qp}{z}$& $p, -p$ & $1$ & $z^2 +p - p^2$ & $p, 1-p$ & $1$ \\ 
\hline
$\frac{q-p}{p} z -\frac{qp}{z}$ & $p, -p$ & $2$ & $z^2 +p - p^2$ & $p, 1-p$ & $1$\\ 
\hline
$\frac{2q}{q^2 -1} z - \frac{p^2 (q^2 +1)}{q (q^2 -1) z}$ & $p, p/q, -p, -p/q$ & $4$ & $z^2 +p - p^2$ & $p, 1-p$ & $1$ \\ 
\hline
$\frac{q+p}{p} z -\frac{qp}{z}$& $p, -p$ & $1$ & $z^2 - (p^2 + p +1)$ & $p, -p-1$ & $2$ \\ 
\hline
$\frac{q-p}{p} z -\frac{qp}{z}$& $p, -p$ & $2$ & $z^2 - (p^2 + p +1)$ & $p, -p-1$ & $2$ \\ 
\hline
$\frac{2q}{q^2 -1} z - \frac{p^2 (q^2 +1)}{q (q^2 -1) z}$ & $p, p/q, -p, -p/q$ & $4$ & $z^2 - (p^2 + p +1)$ & $p, -p-1$ & $2$ \\ 
\hline
$\frac{q+p_{\tau}}{p_{\tau}} z -\frac{qp_{\tau}}{z}$ & $p_{\tau}, -p_{\tau}$ & $1$ & $z^2+ c_{\tau} $ & $p_{\tau},p'_{\tau}, p''_{\tau}$ & $3$ \\ 
\hline
$\frac{q-p_{\tau}}{p_{\tau}} z -\frac{qp_{\tau}}{z}$ & $p_{\tau}, -p_{\tau}$ & $2$ & $z^2 + c_{\tau} $ & $p_{\tau},p'_{\tau}, p''_{\tau}$ & $3$ \\ 
\hline
$\frac{2q}{q^2 -1} z - \frac{p_{\tau}^2 (q^2 +1)}{q (q^2 -1) z}$ & $p_{\tau}, p_{\tau}/q, -p_{\tau}, -p_{\tau}/q$ & $4$ & $z^2 + c_{\tau}$ &$p_{\tau},p'_{\tau}, p''_{\tau}$& $3$ \\ 
\hline
\end{tabular}
\end{center}
PL represents the period length of the orbit of the periodic point with $c_{\tau}=- \frac{\tau^6 + 2 \tau^5 + 4 \tau^4 + 8 \tau^3 + 9 \tau^2 + 4 \tau +1}{4 \tau^2 (\tau +1)^2}$, 
$p_{\tau}= \frac{\tau^3 + 2 \tau^2 + \tau + 1}{2 \tau (\tau +1)}$, $ p'_{\tau}= \frac{\tau^3 - \tau - 1}{2 \tau (\tau +1)}$, and $ p''_{\tau} =-\frac{\tau^3 + 2 \tau^2 + 3 \tau + 1}{2 \tau (\tau +1)} \textrm{ with } \tau \in \Q \setminus \{-1, 0 \}.$

\end{Theorem}
A natural question to pose now is: can we find the triples $(k, b, c)$ such that there exists a rational periodic point $p$ for both $f(z)$ and $\phi_{k, b}(z)$ with 
$$ | \Orb_{f} (p) \cap \Orb_{ \phi_{k, b}} (p) | \ge 2?$$

In fact, examining the table of Theorem \ref{Thm1} gives the following result. 
\begin{Theorem}
Assuming Conjecture \ref{Con:Poonen} and Conjecture \ref{Con:Manes}, if there exists a rational periodic point $p$ such that    
$$ | \Orb_{f} (p) \cap \Orb_{ \phi_{k, b}} (p) | \geqslant 2,$$
then the triples $(k, b, c)$, $c\in\Q, k,b\in\Q^{\times}$, are given as follows:
$$ (k, b, c)= \left(\pm \frac{2 p (p+1)}{2p+1}, \mp\frac{p (p+1) (p^2 + (p+1)^2)}{2p+1}, - (p^2 + p +1) \right) $$
where $p \in \mathbb{Q} \setminus \{ 0, -1/2, -1\}$, and $\Orb_{f} (p) \cap \Orb_{ \phi_{k, b}} (p)= \{p, -p-1 \}$; or

$$ (k, b, c)=\left( \frac{2m_{\tau}}{m_{\tau}^2 -1}, \text{  } -x_i^2  \frac{(m_{\tau}^2 + 1)}{m_{\tau} (m_{\tau}^2 -1)}, - \frac{\tau^6 + 2 \tau^5 + 4 \tau^4 + 8 \tau^3 + 9 \tau^2 + 4 \tau +1}{4 \tau^2 (\tau +1)^2}\right) $$
where $m_{\tau}=\pm x_i/x_j$, $i< j$, and $\Orb_{f} (x_i) \cap \Orb_{ \phi_{k, b}} (x_i)= \{x_i, x_j \} $ where $x_i$, $i=1,2,3$, is defined as in Proposition \ref{prop:1} and
 $\tau \in \mathbb{Q}\setminus\{-1,0\}$. 
 
 In particular, 
$|\Orb_{f} (p) \cap \Orb_{ \phi_{k, b}} (p)|$ cannot exceed $2$ for any periodic point $p\in\Q$.
\end{Theorem}

\begin{Example}
The maps
$f(z)=z^2 -13$ and $  \phi(z)=\frac{24 z}{7} - \frac{300}{7 z}$
 have the following cycles
$ (3, -4)$  and $ (3, -4, -3, 4)$,  respectively.

 The maps $f'(z)=z^2 -\frac{301}{144}$  and $ \phi'(z)= - \frac{115 z}{252} + \frac{31855}{36288 z}$
 have the following cycles
$ \left( \frac{19}{12}, \frac{5}{12}, - \frac{23}{12}\right)$  and $ \left( \frac{5}{12}, \frac{23}{12}, -\frac{5}{12},  - \frac{23}{12} \right),$ respectively.
\end{Example}

\section{Simultaneous Rational Periodic Points of  $\phi_{k_1,b_1} \text{ and } \phi_{k_2,b_2}$}
The aim of this section is to characterize the tuples $(k_1, b_1, k_2, b_2)$ such that the maps $\phi_{k_1,b_1}(z)=k_1z +\frac{b_1}{z}$ and $ \phi_{k_2,b_2}(z)=k_2z +\frac{b_2}{z}$ with  $ k_1, b_1, k_2, b_2 \in \mathbb{Q}^{\times}$ share a common rational periodic point.
Furthermore, we will find the tuples $(k_1, b_1, k_2, b_2)$ such that   
$ | \Orb_{\phi_{k_1, b_1}} (p) \cap \Orb_{ \phi_{k_2, b_2} } (p) | \ge 2$ for some rational periodic point $p$.

A careful work out of the possibilities in Proposition \ref{prop:2} implies the following result.
\begin{Proposition}
\label{prop}
Assuming Conjecture \ref{Con:Manes},  the tuples $(k_1, b_1, k_2, b_2)$ such that there is $p\in\Q^{\times}$ with $p$ a periodic point for the maps $\phi_{k_1,b_1}(z)$ and $\phi_{k_2,b_2}(z)$ are given as follows 

\begin{center}
\begin{tabular}{ | m{7.2em} | m{2.3cm}| m{0.5cm}| m{3cm}| m{2.3cm}| m{0.5cm}|  } 
\hline
$(k_1, b_1)$ & Periodic Points &PL  &$(k_2, b_2)$ & Periodic Points & PL \\ 
\hline
$(1-s_1, s_1 \cdot p^2 )$& $p, -p$ & $1$ & $(1-s_2, s_2 \cdot p^2 )$ & $p, -p$ & $1$ \\ 
\hline
$(s_1 -1, - s_1 \cdot p^2 )$& $p, -p$ & $2$ & $(s_2 -1, - s_2 \cdot p^2 )$ & $p, -p$ & $2$ \\ 
\hline
$( \frac{2s_1}{s_1^2 -1}, -\frac{ p^2 (s_1^2 +1)}{s_1 (s_1^2 -1)})$ & $p, \frac{p}{s_1}, -p, -\frac{p}{s_1}$ & $4$ & $( \frac{2s_2}{s_2^2 -1}, -\frac{ p^2 (s_2^2 +1)}{s_2 (s_2^2 -1)})$ & $p, \frac{p}{s_2}, -p, -\frac{p}{s_2}$ & $4$ \\ 
\hline
$(1-s_1, s_1 \cdot p^2 )$& $p, -p$ & $1$ & $(s_2 -1, - s_2 \cdot p^2 )$ & $p, -p$ & $2$ \\ 
\hline
$(1-s_1, s_1 \cdot p^2 )$& $p, -p$ & $1$ & $( \frac{2s_2}{s_2^2 -1}, -\frac{ p^2 (s_2^2 +1)}{s_2 (s_2^2 -1)})$ & $p, \frac{p}{s_2}, -p, -\frac{p}{s_2}$ & $4$ \\ 
\hline
$(s_1 -1, - s_1 \cdot p^2 )$& $p, -p$ & $2$ & $( \frac{2s_2}{s_2^2 -1}, -\frac{ p^2 (s_2^2 +1)}{s_2 (s_2^2 -1)})$ & $p, \frac{p}{s_2}, -p, -\frac{p}{s_2}$ & $4$ \\ 
\hline
\end{tabular}
\end{center}
where PL is the period length of the orbit of the periodic point, $s_i\in\Q\setminus\{0,1\}$ if PL is either $1$ or $2$; and $s_i\in\Q\setminus\{0,\pm1\}$ if PL is $4$.
\end{Proposition}

Using Proposition \ref{prop}, one obtains the following result.
\begin{Theorem}
\label{thm:same}
Assuming Conjecture \ref{Con:Manes}, the quadruples $(k_1,b_1,k_2,b_2)$ such that there exists a periodic point $p\in\Q^{\times}$ and maps
$\phi_{k_i,b_i}(z)=k_iz +\frac{b_i}{z}$, 
with $k_i, b_i\in \mathbb{Q}^{\times}, i=1,2,$ satisfying
$$ | \Orb_{\phi_{k_1, b_1}} (p) \cap \Orb_{ \phi_{k_2, b_2} } (p) | \ge 2$$ are given as follows:

$$ (k_1, b_1, k_2, b_2)= (s_1 -1, -s_1 \cdot p^2, s_2 -1, -s_2 \cdot p^2 ), \qquad s_1, s_2 \in \mathbb{Q} \setminus \{0, 1 \};\textrm{ or}$$
$$ (k_1, b_1, k_2, b_2)= \left(s_1-1 , -s_1 \cdot p^2, \frac{2s_2}{(s_2^2 -1)}, - p^2 \cdot \frac{(s_2^2 + 1)}{s_2 (s_2^2 -1)}\right), \qquad s_1 \in \mathbb{Q} \setminus \{0, 1 \},\, s_2 \in \mathbb{Q} \setminus \{ 0, \pm1\},\,\textrm{ or }$$
$$ (k_1, b_1, k_2, b_2)= \left( \frac{2s_1}{s_1^2 -1} , - p^2 \cdot \frac{(s_1^2 + 1)}{s_1 (s_1^2 -1)}, \frac{2s_2}{s_2^2 -1}, - p^2 \cdot \frac{(s_2^2 + 1)}{s_2 (s_2^2 -1)} \right), \qquad
s_1 \neq \pm s_2,\textbf{ } s_1, s_2 \in \mathbb{Q} \setminus \{0, \pm1\} $$
where in these cases $\Orb_{\phi_{k_1, b_1}}(p) \cap \Orb_{\phi_{k_2, b_2}}(p) = \{p, -p\}.$

In particular, if $ | \Orb_{\phi_{k_1, b_1}} (p) \cap \Orb_{ \phi_{k_2, b_2} } (p) | = 4$, then $\phi_{k_1, b_1}= \pm \phi_{k_2, b_2}$.
\end{Theorem}
\begin{Example}
The maps $\phi(z)=\frac{4z}{3} - \frac{3}{10 z}$ and $ \phi(z)=-\frac{3z}{4} + \frac{27}{20 z}$ have the following cycles $\left(\frac{3}{5}, \frac{3}{10}, -\frac{3}{5}, -\frac{3}{10}\right)$ and , $ \left(\frac{3}{5}, \frac{9}{5}, -\frac{3}{5},- \frac{9}{5}\right)$, respectively. 
\end{Example}
We now recall the main theorem in \cite{BakerDeMarco}.
\begin{Theorem}
Let $d \geqslant 2$ be an integer. Fix $c_1, c_2 \in \mathbb{C}$. The set of $t \in \mathbb{C}$ such that both $c_1$ and $c_2$ are preperiodic for $z^d + t$ is infinite if and only if $c_1 ^d = c_2 ^d$. 
\end{Theorem}
Assuming Conjecture \ref{Con:Manes}, we will prove the following theorem in this section. 
\begin{Theorem} \label{the last theorem}
Assume Conjecture \ref{Con:Manes} holds. 
Let $a, b \in \mathbb{Q}^{\times}$. 
There exists infinitely many rational pairs $(t_1, t_2)\in \mathbb{Q}^{\times} \times \mathbb{Q}^{\times}$ such that $a$ and $b$ are both rational periodic points of $\phi_{t_1, t_2}(z)=t_1 \cdot z + \frac{t_2}{z}$ if and only if $a^2 = b^2$.
\end{Theorem}
Before we start proving Theorem \ref{the last theorem}, we give the following straightforward application of Proposition \ref{prop:2}.

\begin{Lemma} \label{periodiconetwo}
If $\phi_{k, b}(z)=kz+\frac{b}{z}$ has two distinct rational periodic points $q_1, q_2\in\Q^{\times}$ with exact periods $1$ and $2$, respectively, where $k, b \in \mathbb{Q}^{\times}$, then 
$$ \phi_{k, b}(z)= \frac{- q_2^2 -q_1^2}{ q_2^2 -q_1^2}z + \frac{2q_1^2 q_2^2}{(q_2^2 -q_1^2) z}.$$
\end{Lemma}

The following Lemma is \cite[Proposition 10]{Manes}. 

\begin{Lemma} \label{periodicwith4}
Let $\phi_{k, b}(z)=kz+\frac{b}{z}$ where $k\in\Q\setminus\{0,-1/2\}, b \in \mathbb{Q}^{\times}/(\mathbb{Q}^{\times})^2$. Then
 $\phi_{k, b}(z)$ cannot have two periodic points with exact period $n$ and exact period $4$,  $n=1,2$.
\end{Lemma}
 
 Now we proceed with the proof of Theorem \ref{the last theorem}
\begin{ProofOf}{Theorem \ref{the last theorem}}
Suppose that there exists infinitely many rational pairs $(t_1, t_2)$ such that $a$ and $b$ are both rational periodic points of $\phi_{t_1, t_2}(z)$. By Lemma \ref{periodiconetwo} and Lemma \ref{periodicwith4}, $a$ and $b$ cannot have orbits with different period lengths under the action of $\phi_{t_1, t_2}(z)$. We may therefore assume that both $a$ and $b$ have period $n$,  where $n$ is either $1, 2$ or $4$. In view of Proposition \ref{prop}, we must have $a= -b$. 
\end{ProofOf}
 
 \section{Simultaneous rational periodic points of two quadratic polynomials}
 
 One remarks that a rational number $p$ can be a periodic point for infinitely many non-linearly conjugate maps of the form $\phi_{k,b}(z)=kz+b/z$, $k,b\in\Q^{\times}$, see Proposition \ref{prop}. However, this is not the case for quadratic polynomial maps.
\begin{Theorem} \label{Thmatmost3}
Assume Conjecture \ref{Con:Poonen} holds. If $q\in\Q$ is a periodic point for every map $f_i(z)=z^2+c_i$, $i=1,\cdots,m$, $c_i \in \Q$, $c_i \neq c_j \textrm{ if } i \neq j$, then $m\le 3$. 
\end{Theorem}
\begin{Proof}
Suppose $q$ is rational fixed point of $f(z)=z^2 +c$. It follows that $c$ is $\frac{1}{4} -\rho^2$ for some $\rho \in \mathbb{Q}$ where $q$ is either $1/2 + \rho$ or $1/2 - \rho$, see Proposition \ref{prop:1},  so $\rho=\pm(q - 1/2)$, and $c$ is determined uniquely. 

Similarly, if $q$ is a rational periodic point of exact period $2$ of $f(z)=z^2 +c$, then $c$ is $-\frac{3}{4} -\sigma^2$ for some $\sigma \in \mathbb{Q}^{\times}$. In this case, $q$ is either $-1/2 + \sigma$ or $-1/2 - \sigma$, see Proposition \ref{prop:1}, so that $\sigma=\pm(q + 1/2)$, and $c$ is determined uniquely.

We now suppose $q$ is rational periodic point of exact period $3$ of $f(z)=z^2 +c$. Then, 
$c=- \frac{\tau^6 + 2 \tau^5 + 4 \tau^4 + 8 \tau^3 + 9 \tau^2 + 4 \tau +1}{4 \tau^2 (\tau +1)^2},$
for some $\tau \in \mathbb{Q} \setminus \{ -1, 0\}$, and $q$ is one of the $x_i$'s in Proposition \ref{prop:1}.  

Assume that $q=x_1(\tau)= \frac{\tau^3 + 2 \tau^2 + \tau + 1}{2 \tau (\tau + 1)}.$ This implies that $\tau$ is a root of the following cubic polynomial equation  
\begin{equation*}
x^3 + (2 -2q) x^2 + (1 -2q) x + 1=0.
\end{equation*}
We now check whether the latter polynomial has a rational root other than $\tau$. We know that for other roots $x_1, x_2$ of this polynomial we have 
$$x_1 + x_2= 2q - 2 - \tau, \textbf{ } x_1 \cdot x_2 = -1/\tau.$$ 
Therefore $x_1$ and $x_2$ satisfy the following quadratic polynomial equation
$$ x^2 + (\tau + 2 - 2q) x - \frac{1}{\tau}=0.$$
This polynomial has rational roots if and only if the discriminant
$ ( \tau + 2 - 2q)^2 + \frac{4}{\tau}$
is a square in $\mathbb{Q}$. Writing $q$ in terms of $\tau$, the latter condition on the discriminant is equivalent to the fact that the following expression 
$ \frac{\tau^4 + 6 \tau^3 + 7 \tau^2 + 2 \tau + 1}{\tau^2 (\tau +1)^2} $
is a square in $\mathbb{Q}$. Equivalently, one needs to find a rational solution on the elliptic curve 
\begin{equation}\label{curve}y^2=\tau^4 + 6 \tau^3 + 7 \tau^2 + 2 \tau + 1.\end{equation}
Using \Magma, \cite{Magma}, one can check $(0, 1), (0, -1), (-1, 1), (-1, -1)$, and the two points at infinity are the only rational points on this elliptic curve. Since $\tau \neq 0, -1$, the discriminant cannot be a square in $\mathbb{Q}$. 

Arguing similarly, if $q=x_2(\tau)= \frac{\tau^3 - \tau - 1}{2 \tau (\tau +1)},$ then $\tau$ 
 is a root of the following polynomial
$
x^3 -2q x^2 - (1 +2q) x - 1$.
The other two roots $x_1$ and $x_2$ satisfy the quadratic polynomial
$ x^2 + (\tau - 2q) x + \frac{1}{\tau}.$
This polynomial has rational roots if and only if the discriminant
$ ( \tau - 2q)^2 - \frac{4}{\tau}$
is a square in $\mathbb{Q}$. Writing $q$ in terms of $\tau$, the latter occurs if and only if
$ \frac{\tau^4 - 2 \tau^3 - 5 \tau^2 - 2 \tau + 1}{\tau^2 (\tau +1)^2} $
is a square in $\mathbb{Q}$. Knowing that the rational points on the elliptic curve defined by 
$y^2=\tau^4 - 2 \tau^3 -5 \tau^2 - 2 \tau + 1$
are $(0, 1), (0, -1), (-1, 1), (-1, -1)$, and the two points at infinity, \Magma \cite{Magma}, the discriminant cannot be a square in $\mathbb{Q}$. 

Finally, assuming
$q=x_3(\tau)= - \frac{\tau^3 + 2 \tau^2 + 3 \tau + 1}{2 \tau (\tau + 1)}$, for some $\tau\in\Q\setminus\{-1,0\}$, leads us to studying the rational points on the elliptic curve 
$y^2=\tau^4 + 2 \tau^3 + 7 \tau^2 + 6 \tau + 1$ which is $\Q$-birational to the elliptic curve in (\ref{curve}). Hence, we conclude that if $q= x_i(\tau)$, $i=1,2,3$, for some $\tau\in \Q \setminus \{-1, 0\}$, then $\tau$ is uniquely determined.  This fact together with the fact that $x_1$, $x_2$ and $x_3$ are permuted under the action of the quadratic map imply that if $q$ is a rational periodic point of exact period $3$ of $f(z)=z^2 +c$, then both the orbit of $q$ and $c$ are determined uniquely. 
\end{Proof}

\begin{Remark}
In the proof of Theorem \ref{Thmatmost3}, we showed that for distinct $c_1, c_2 \in \Q$, given a rational point $p$ such that $p$ is a periodic point of $z^2+c_1$ and $z^2+c_2$ of exact period $m\le 3$ and $n\le 3$, respectively, then $m \ne n$. 

A rational point $p \in \Q$ is a  periodic point of the maps $f(z)= z^2 + p -p^2$ and $f'(z)= z^2 - (p^2 + p + 1)$ of exact periods $1$ and $2$, respectively. 

If $p=p_{\tau} \in \Q$ is a periodic point of exact period $3$ of the map $f_3(z)=z^2 + c_{\tau}$, then $p_{\tau}$ is a periodic point of exact period $1$ of $f_1(z)= z^2 + p_{\tau} -p_{\tau}^2$, and a periodic point of exact period $2$ of $f_2(z)= z^2 - (p_{\tau}^2 + p_{\tau} + 1)$.

\end{Remark}

\begin{Example}
The maps $f_1(z)=z^2 - \frac{6161}{1600}$, $f_2(z)=z^2 - \frac{15841}{1600}$, and $f_3(z)= z^2 - \frac{7841}{1600}$ have the following cycles $\left(\frac{101}{40}\right)$, $ \left(\frac{101}{40}, - \frac{141}{40} \right)$, and $ \left(\frac{101}{40}, \frac{59}{40}, - \frac{109}{40} \right)$, respectively. 
\end{Example}

\end{document}